\documentclass[12pt]{amsart}
\usepackage{amsfonts, amssymb, amsmath, latexsym}

\newcommand{\la}{\langle}
\newcommand{\ra}{\rangle}

\newcommand{\goth}{\mathfrak}

\newcommand{\pz}{\partial}
\newcommand{\C}{\mathbb{C}}
\newcommand{\R}{\mathbb{R}}
\newcommand{\bz}{\bar{\partial}}
\newcommand{\px}{\frac{\partial}{\partial x}}
\newcommand{\py}{\frac{\partial}{\partial y}}
\newcommand{\Pz}{\frac{\partial}{\partial z}}
\newcommand{\Bz}{\frac{\partial}{\partial \bar{z}}}
\newcommand{\bb}{\mathbb}
\newcommand{\Fl}{\mathcal{F}} 
\newtheorem{theorem}{Theorem}[section]
\newtheorem{proposition}{Proposition}[section]
\newtheorem{corollary}[theorem]{Corollary}
\newtheorem{lemma}[theorem]{Lemma}
\newtheorem{remark}{Remark}[section]

\begin{document}

\title{ surfaces in $\bb{S}^4$ with normal harmonic Gauss map } 
\author{Eduardo Hulett}
 \thanks{Partially supported by ANPCyT, CONICET and SECYT-UNC,
Argentina.} 
\address{C.I.E.M. - Fa.M.A.F.
 Universidad Nacional de C\'ordoba,
Ciudad Universitaria,
5000 C\'ordoba, Argentina.
-phone/fax: +54 351 4334052/51}
\email{ hulett@famaf.unc.edu.ar}
\date{}

 \begin{abstract}
 We consider  conformal immersions of Riemann surfaces   in $\bb{S}^4$ and study their Gauss maps with values in the Grassmann  bundle $\Fl = SO_5/T^2 \to \bb{S}^4$. The energy of maps from Riemann surfaces into $\Fl$  is considered with respect to   the normal metric  on the target and immersions  with harmonic Gauss maps are characterized. We also show that the normal-harmonic map equation for Gauss maps is a completely integrable system, thus giving  a partial answer of a question posed by Y. Ohnita in~\cite{ohnita}. Associated $\bb{S}^1$-families  of parallel mean curvature immersions in $\bb{S}^4$  are considered. A  lower bound of the  normal energy of Gauss maps is obtained in terms of the genus of the surface.  
 \end{abstract}
\maketitle

\section{Introduction }
 A conformally immersed Riemann surface $f: \Sigma \to \bb{S}^4$ has a well defined Gauss map which asigns to a point $ p \in \Sigma $ the oriented $2$-plane $df (T_p \Sigma)$ viewed as  a point  in the Grassmann bundle  $G_2(T \bb{S}^4)$ of oriented $2$-planes in $T\bb{S}^4$, which is a       homogeneous manifold diffeomorphic to the flag manifold $\Fl = SO_5 /SO_2 \times SO_2$.    
An early result of Eells and Salamon~\cite{eells-salamon} asserts that a smooth map $\phi :\Sigma \to \Fl$ of a Riemann surface is {\it primitive}  if and only if it is the Gauss map of a minimal (i.e. conformal harmonic) immersion  $\pi \circ \phi :\Sigma \to \bb{S}^4$, where $\pi : \Fl \to \bb{S}^4$ is   the homogeneous projection. 
  Primitive maps from Riemann surfaces in $\Fl$ are  holomorphic maps with respect to the horizontal $F$-structure  determined by the canonical $4$th-order  Cartan inner authomorphism of the Lie algebra $\goth{sl}_5(\C)$~\cite{burstall-pedit}. A remarkable property  of  primitive maps is   their {\it equi-harmonicity} i.e.  harmonic with respect to every  invariant metric on the target   $\Fl$~\cite{black},~\cite{bolton-woodward},~\cite{burstall-pedit}. \\
  A related question is  to study the geometry of immersions of orientable  surfaces in $\bb{S}^4$ for which their  Gauss maps satisfy the harmonic map equation with respect to a specific $SO_5$-invariant metric on $\Fl$. 
In this article we  consider  conformal  immersions $f:\Sigma \to \bb{S}^4$ of Riemann surfaces 
whose  Gauss maps are {\it normal-harmonic} or harmonic with respect to  the so-called normal metric on the target $\Fl$.    Since naturally reductive homogeneous spaces may be considered as generalizations of Riemannian symmetric spaces, the choice of the normal metric on $\Fl$ seems to be very natural.  There are however   other geometrically interesting  possible choices of invariant metrics on $\Fl$  which are considered in~\cite{hulettinprogress}. \\
The  first main result of the paper is Theorem~\ref{harmonicgausslift} which is  a generalization of the well-known Theorem of Ruh-Vilms~\cite{ruh-vilms}. It characterizes immersed surfaces in $\bb{S}^4$ with normal-harmonic Gauss maps.
Although this result is  quite natural and elementary,  to the best of this
author's knowledge it has not been reported before. \\
Theorem~\ref{comp.int.Gaussmap}, the second main result in the paper,  establishes the complete integrability of the normal-harmonic map equation for Gauss maps. 
In other words it asserts  that   the normal-harmonic map equation for Gauss maps $\Sigma \to \Fl$ can be encoded in a loop of flat connections $\bb{S}^1 \ni \lambda \mapsto d+ \alpha_{\lambda}$  on the trivial principal $SO_5$ bundle over $\Sigma$, a manifestation of complete integrability.
Its proof is consequence of  identity~\eqref{specialproperty} in Lemma~\ref{specialpropertylemma}, which is a special property of the Gauss map. This fact  provides a partial  answer to the following  question posed by Y. Ohnita~\cite{ohnita}:     Are there examples of harmonic maps, other than  super-horizontal and primitive maps  from a Riemann surface $\Sigma$  to a $k$-symmetric ($k >2$)  manifold  $G/K$ which satisfy condition of Lemma~\ref{specialpropertylemma}?    \\
     
     The paper is organized as follows. In section~\ref{structure}  we derive the structure equations of isometric immersions of orientable surfaces in $\bb{S}^4$. As is well known, any such immersion determines a conformal or Riemann surface structure on an orientable surface for which the immersion results conformal.      
     In section~\ref{gaussmap} we give  a detailed construction of the Gauss map and some elementary facts on the geometry of the flag manifold $\Fl \equiv SO_5/SO_2 \times SO_2$. We describe also the Maurer-Cartan one form or moment map $\beta$ of $\Fl$ following~\cite{burstall-rawnsley}. 
     In section~\ref{normalharmonicmaps} we  closely follow~\cite{black} and~\cite{burstall-rawnsley} to  derive  the harmonic map equation for smooth maps into $\Fl$ when  the normal metric has been fixed,  we call these maps normal-harmonic. When maps from Riemann surfaces are considered, the normal-harmonic map equation takes a very simple form.
     In section~\ref{tensionGauss} we derive an explicit formula  of the tension of the Gauss map of a conformal immersion in terms of the mean curvature vector of the immersion. A direct  consequence is the proof of   Theorem~\ref{harmonicgausslift}.   Also we obtain Lemma~\ref{specialpropertylemma} which establishes property~\eqref{specialproperty},  a distinctive algebraic-geometric property of Gauss maps. As a consequence of  we show that  the normal-harmonic map equation for Gauss maps is a completely integrable system. We use this information to describe the $\bb{S}^1$-loop of conformal immersions of a simply connected Riemann surface determined by a normal-harmonic Gauss map.  
      The last section deals with a  computation of the normal energy of Gauss maps. We  obtain a formula relating  the normal energy of the Gauss map and  the Willmore energy of the corresponding immersion. This allows  to  obtain  a lower bound for the normal energy of Gauss maps which depends on the genus of the immersed Riemann surface.

\section{Structure equations }\label{structure} 

On  $\R^5$    with coordinates $(x_1, x_2, x_3, x_4, x_5)$   the euclidean  metric 
 $$\la .,.\ra = dx^2_1 + dx^2_2 + dx^2_3 + dx^2_4 + dx^2_5$$
induces the usual  canonical Riemannian metric $\langle .,. \rangle $ with constant curvature one on 
the  unit sphere   $
      \bb{S}^4 = \{ x \in  \R^5 : \langle x,x \rangle =1 \}.$ 
The  matrix  Lie group $SO_5= \{A \in Gl_5 (\R) : A^t A = I, det A =1 \}$ acts transitively on  $\bb{S}^4$ by isometries.   \\ 
 An immersion $f : \Sigma \to \bb{S}^4$  of a Riemann surface is   conformal if  
   $\la 
  f_z ,  f_z \ra^{c}=0$, for every  local complex
  coordinate $z=x+iy$ on $M$, where 
  $$\Pz = \frac{1}{2} ( \frac{\partial}{\partial x}- i \frac{\partial}{\partial y}), \quad \Bz = \frac{1}{2} ( \frac{\partial}{\partial x} + i \frac{\partial}{\partial y}),$$
are the  complex partial derivatives  and  $\la \, , \, \ra^{c}$ is the complex bilinear  extension of
   the riemannian  metric to $\C^5$:
   $$  \la z,w \ra^c =  z_1 w_1 + z_2 w_2 + z_3 w_3 + z_4 w_4+z_5 w_5.  $$
   
    Thus $ f$   is conformal if and only if    on any   local complex coordinate $z =x + iy$ it satisfies 
   \begin{equation} \label{confID}
       \la f_x , f_y \ra =0, \quad \| f_x \|^2 =\| f_y \|^2.
  \end{equation}
 
 We fix on $\Sigma$ the induced Riemannian metric $g= f^* \la \,,\, \ra$ so that  $ f: (M,g) \to \bb{S}^4$ is an isometric immersion. 
    The {\it $2$nd Fundamental form} of the  surface $f:M \to \bb{S}^4$  is defined by   $II =  d^2f$ i.e. 
    $$ II  = - \la df, d N_1 \ra N_1 -  \la df, d N_2 \ra N_2, $$
    where  $\{ N_1, N_2 \}$ is a local  orthonormal  frame.   
    On any local  chart $(U, z=x+iy)$ of $M$ we introduce   
          a conformal parameter $u$ defined by $\la  f_z,  f_z \ra=e^{2u}$, so that    $g|_U=2 e^{2u}( dx^2 +   dy^2)$. 
 The {\it mean curvature} vector of $f$ is  defined by $H=  \frac{1}{2}trace II $, which   in terms of $f$ and $u$ is   given by 
      $$ H= e^{-2u}  f^{\bot}_{\bar{z} z}. $$
Defining  $h_i = \la H, N_{i} \ra = e^{-2u}  \la f_{\bar{z} z}, N_{i} \ra $ we decompose  $     H = h_1 N_1 + h_2 N_2$. \\
Since $f$ is conformal 
   $$ 
    2\la  f_{\bar{z} z} ,  f_z \ra^{c} = \Pz \la  f_z,  f_z
    \ra^{c} =0, \quad 
    2\la  f_{\bar{z} z} , f_{\bar{z} } \ra^{c} = \Pz \la f_{\bar{z} }, f_{\bar{z} }
    \ra^{c} =0,
     $$
    thus      $f_{\bar{z} z}$ has no tangential component and so   
    \begin{equation}\label{Laplacef}   f_{\bar{z} z} = -e^{2u} f + e^{2u}H.  
    \end{equation}
    
    Defining $\xi_i : = \la f_{zz}, N_{i} \ra^c$,  we decompose  
    $$  f_{z z} = 2  u_z .  f_z + \xi_1 N_1 + \xi_2 N_2.  $$
From  this equation we obtain
\begin{equation}\label{holfundform} 
    f_{z z}^{\bot} =  \frac{1}{4} \left[ II( \px, \px ) - II (\py, \py) \right] - \frac{i}{2} II( \px, \py) = II (\Pz, \Pz ).   
    \end{equation}

\vskip .5cm  

 On the other hand from 
     \begin{equation}
     \begin{array}{ll}
      0= \Pz \la N_1, f_z \ra = \la \Pz N_1, f_z \ra + \la N_1, f_{\bar{z}z} \ra,\\
      \\
      0= \Pz \la N_1, f_{\bar{z}} \ra = \la \Pz N_1, f_{\bar{z}} \ra + \la N_1, f_{\bar{z}\bar{z}} \ra,\\
      \end{array}
      \end{equation}
     we obtain 
     \begin{equation}
     \begin{array}{ll}
   \la \Pz N_1, f_z \ra = - e^{2u} h_1. \\
   \\
   \la \Pz N_1, f_{\bar{z}} \ra= -\la N_1, f_{zz} \ra^c = - \xi_1 N_1. 
   \end{array}
   \end{equation}
In an analogous way
\begin{equation}
\begin{array}{ll}         
\la \Pz N_2, f_z \ra = -e^{2u} h_2, \\
\\
\la \Pz N_2, f_{\bar{z}} \ra= -\la N_2, f_{zz} \ra^c = - \xi_2 N_2.
\end{array}
\end{equation}   
Defining $\sigma:= \la \Pz N_2 , N_1 \ra$, we obtain 
 the complex derivative of the normal fields 
 \begin{equation}
 \begin{array}{ll}
  \Pz N_1 = -h_1 f_z - e^{-2u} \xi_1.
    f_{\bar{z}} - \sigma N_2.\\
    \\
   \Pz N_2 = - h_2 f_z - e^{-2u} \xi_2.  f_{\bar{z}} + \sigma  N_1
    \end{array} 
 \end{equation}
    These equations together are the structure equations of the conformal immersion $f:\Sigma \to \bb{S}^4$, 
     \begin{equation}\label{structureeqs}
      \left \{ 
      \begin{array}{ll} 
         f_{\bar{z} z} = -e^{2u} f + e^{2u}H \\
        f_{z z} = 2  u_z .  f_z + \xi_1 N_1 + \xi_2 N_2, \\
        \,  \Pz N_1 = -h_1 f_z - e^{-2u} \xi_1.
    f_{\bar{z}} - \sigma N_2.\\
    \Pz N_2 = - h_2 f_z - e^{-2u} \xi_2.  f_{\bar{z}} + \sigma  N_1.   \\       
       \end{array} \right.
       \end{equation}
       
       The integrability condition of these equations are  the Gauss, Codazzi and Ricci  equations respectively,
       \begin{equation}\label{gausscodazziricci}
     \left \{ 
       \begin{array}{ll}
       \text{(G)} & 2 u_{\bar{z}z} = e^{-2u} (|\xi_1|^2 + |\xi_2|^2) - e^{2u}( 1 + \| H\|^2),\\
        \text{(C1) } & e^{2u} ( \Pz h_1 + h_2 \sigma) = \Bz \xi_1 + \xi_2 \bar{\sigma},\\
      \text{(C2) } & e^{2u} ( \Pz h_2 - h_1 \sigma) = \Bz \xi_2 - \xi_1 \bar{\sigma},\\
 \text{(R)}&      - Im (  \sigma_{\bar{z}} ) = e^{-2u} Im (\xi_1 \bar{\xi}_2).\\
      \end{array} \right.
       \end{equation}
   
Let $\nabla^{\bot}$ be the connection on the normal bundle $\nu$ of the immersed surface. Then from the compatibility equations we obtain
\begin{equation}\label{normalderivH}
  \nabla^{\bot}_{\Pz} H =( \Pz h_1 + h_2 \sigma) N_1 +   ( \Pz h_2 - h_1 \sigma) N_2.
\end{equation}

The Codazzi equations  $(C1)$ and $(C2)$ in~\eqref{gausscodazziricci} may be expressed also in the following form  
\begin{equation}\label{covderivH}
 \nabla^{\bot}_{\Pz} H = e^{-2u} (\Bz \xi_1 + \xi_2 \bar{\sigma}) N_1 + e^{-2u}(\Bz \xi_2 -  \xi_1 \bar{\sigma}) N_2
\end{equation}
\vskip .5cm

The  induced metric 
$g=f^* \la \,,\, \ra$ on $\Sigma$ is given in terms of the conformal parameter $u$  by  $g=2 e^{2u}dz \otimes d\bar{z}$. The  Gaussian curvature of the surface $(\Sigma,g)$ is just the curvature of the induced metric which   is given by 
$$K= - \Delta_g u = -2 e^{-2u}  u_{\bar{z} z}, $$
 where $\Delta_g=  2 e^{-2u} \Bz \Pz$  is the Laplace operator on $\Sigma$ determined by $g$.  Gauss equation (G) in~\eqref{gausscodazziricci} becomes 
 \begin{equation}\label{gausscurvature} 
K = ( 1 + \| H\|^2) - e^{-4u} (|\xi_1|^2 + |\xi_2|^2),
\end{equation}

Let  $\omega = \la \nabla^{\bot} N_2 , N_1 \ra$ be the connection one form of the normal bundle $\nu$.  Then the normal curvature is defined by $d \omega = K^{\bot} dA$, where $dA$ is the area form of the induced metric $g$. Thus in terms of $\sigma$ we obtain $\omega = 2 Re (\sigma dz)$ and so
$$
                 d \omega  =  -4 Im (\sigma_{\bar{z}}) dx \wedge dy $$
Hence since $d \omega = K^{\bot} dA_g$, and $dA_g = 2 e^{2u} dx \wedge dy$, the normal curvature function is given by 
\begin{equation}\label{normalcurvature}
K^{\bot}= -e^{-2u} Im (\sigma_{\bar{z}}).
\end{equation}

Let  $F: \Sigma \to SO_5$ be a (local) adapted frame of $f$ i.e.  
   \begin{equation} \label{defadaptedframe}
      f(x) = F_0(x), \quad df (T_x M)= span \{ F_1(x), F_2(x) \},  x \in M,
      \end{equation}
       where $F_i, 0 \leq i \leq 4$ are the columns of the orthogonal matrix $F$.  By a gauge transformation   (rotating within  the complex line generated by $F_1 - i F_2$) we can assume that 
 \begin{equation}\label{fadapted}
      f_z = \frac{e^{u}}{\sqrt{2}}(F_1 - i F_2).
      \end{equation}
      Using~\eqref{fadapted} and routine  computation we  compute the complex derivative $\Pz F$ of the frame $F$ in terms of the frame itself. 
        
\begin{equation}\label{structureeqs2}
\left \{ 
\begin{array}{l}
\Pz F_0 = \frac{e^{u}}{\sqrt{2}} F_1 - i \frac{e^{u}}{\sqrt{2}} F_2, \\
\\
  \Pz F_1=
-\frac{e^{u}}{\sqrt{2}}.f -i u_z F_2 + ( \frac{ e^{-u} \xi_1+
e^{u}h_1 }{\sqrt{2}}).N_1 +( \frac{ e^{-u} \xi_2+
e^{u}h_2 }{\sqrt{2}}).N_2 
\\
\\
 \Pz F_2= i \frac{e^{u}}{\sqrt{2}}.f + i  u_z F_1 + i  ( \frac{
e^{-u} \xi_1 - e^{u}h_1 }{\sqrt{2}}).N_1 + i  ( \frac{
e^{-u} \xi_2 - e^{u}h_2 }{\sqrt{2}}).N_2  \\
\\
\Pz N_1 = -(\frac{e^{-u} \xi_1+ e^u h_1 }{\sqrt{2}}). F_1 - i (\frac{ e^{-u} \xi_1 - e^u h_1}{\sqrt{2}}).F_2 - \sigma N_2\\
\\
\Pz N_2 = - (\frac{e^{-u} \xi_2+ e^u h_2 }{\sqrt{2}}) F_1 - i (\frac{e^{-u} \xi_ 2 - e^u h_2 }{\sqrt{2}}). F_2 + \sigma N_1 \\
\end{array}
\right.
\end{equation}
\vskip .5cm
Set 
\begin{equation}\label{aibi} 
a_i =  \frac{ e^{-u} \xi_i+
e^{u}h_i }{\sqrt{2}}, \quad  b_i = \frac{  e^{-u} \xi_i - e^u h_i}{\sqrt{2}},\,\, i=1,2
\end{equation}      
System~\eqref{structureeqs2} can be written  in matrix form  as  $  F_z = F A$, where  $F =(F_1, F_2, F_3, F_4, F_5)$ in column notation and
      
\begin{equation}\label{matrA}
A =
\begin{pmatrix}
0  &   -\frac{e^u}{\sqrt{2}}  &   i\frac{e^u}{\sqrt{2}}  &    0  & 0\\

\frac{e^u}{\sqrt{2}}& 0  &  i u_z &  - a_1   & -a_2\\

-i \frac{e^u}{\sqrt{2}} & -i u_z & 0 & -i  b_1  &  -i b_2  \\

0 & a_1  &  i b_1 &       0       &     \sigma        \\

0 & a_2  &  i b_2 &    -\sigma     &        0              \\        
\end{pmatrix},
\end{equation}

\begin{equation}\label{matrB}
B =
\begin{pmatrix}
0  &   -\frac{e^u}{\sqrt{2}}  &   -i\frac{e^u}{\sqrt{2}}  &    0  & 0\\

\frac{e^u}{\sqrt{2}}& 0  &  -i u_{\bar{z}} &  - \bar{a}_1   & -\bar{a}_2\\

i \frac{e^u}{\sqrt{2}} & i u_{\bar{z}} & 0 & i  \bar{b}_1  &  i \bar{b}_2  \\

0 & \bar{a}_1  &  -i \bar{b}_1 &       0       &     \bar{\sigma}        \\

0 & \bar{a}_2  &  -i \bar{b}_2 &    -\bar{\sigma}     &        0              \\       
\end{pmatrix}.
\end{equation}       
       
 Since $\bar{A} =B$ we get $ F_{\bar{z}} = F B$, so that in  terms of matrices $A,B$ the integrability condition $ F_{z \bar{z}} =  F_{\bar{z}z}$, 
  is equivalent to
  \begin{equation}\label{equationAB}
        A_{\bar{z}} -  B_z = [A,B], 
  \end{equation}
 which encodes  the equations of Gauss, Codazzi and Ricci of the immersion   given before. 
 Let  $\Theta$ denote the left Maurer-Cartan form of the group $SO_5$, and consider   the   pullback  $\alpha = F^* \Theta = F^{-1} dF = A dz + B d\bar{z}$. Then equation~\eqref{equationAB} is equivalent to   
 \begin{equation}
                   d \alpha + \frac{1}{2} [ \alpha \wedge \alpha] =0.
 \end{equation}
  
\vskip .5cm

 \section{The Gauss map } \label{gaussmap}

  A  smooth map  $f: \Sigma \to \bb{S}^4 $ determines a map $\varphi : \Sigma \to  \bb{CP}^4$  defined by $\varphi (x) = [f(x)]$.  Let  $L \to \bb{CP}^4$ be  the tautological line bundle whose  fiber over $\ell \in \bb{CP}^5$ is the complex line $\ell$ itself. Then $\varphi$ (and so $f$)  determines a  
complex line subbundle $\ell_0   \subset \Sigma \times \C^5 $ given by 
 \begin{equation}\label{L0}
             \varphi^*(L)=  \ell_0 = \{ (x,v) \in \Sigma \times \C^5 :v \in \C f(x) \} 
 \end{equation}   
 
 Any smooth subbundle $E \subset \Sigma \times \bb{C}^5$  can be equipped with a holomorphic structure for which a local section $s \in \Gamma(E)$ is holomorphic if and only if $\pi_E (\Pz s)=0$, where $\pi_E$ is the orthogonal projection onto $E$ and $z$ is any local complex coordinate on $\Sigma$. Now   $f$ is harmonic if and only if $\varphi$ is harmonic since $\varphi = \iota \circ f$, where $\iota$ denotes the totally geodesic embedding $\bb{S}^4 \hookrightarrow \bb{CP}^4$.     In turn $\varphi$ is harmonic if and only if
  the map 
 $$
                     d \varphi (\Pz): \ell_0 \to \ell^{\bot}_0, \quad s \mapsto \pi^{\bot}_{\ell_0} (\Pz s)  $$
  is holomorphic i.e. it sends holomorphic sections of $\ell_0$ to holomorphic sections of $\ell^{\bot}_0$. 
Since $\la f, f_{z} \ra =0$,  $f$ is a global holomorphic section of $\ell_0$ and $\pi^{\bot}_{\ell_0} (\Pz f) = f_z  $.  Now  if $f$ is conformal,   $f_z$ is a section of $\ell^{\bot}_0 $  satisfying $\pi^{\bot}_{\ell_0} ( \Bz f_z) =      e^{2u} H$ by equation~\eqref{Laplacef}.  In particular  $\varphi$ (and hence $f$) is harmonic iif $H=0$.  \\
Define $\ell_1 := \pi^{\bot}_{\ell_0} (\Pz \ell_0) \subset \ell_0^{\bot} \subset \Sigma \times \C^5$. Then $f$ is conformal  if and only if $\ell_1$ is an isotropic subbundle.    A conformal map $f$ is called isotropic  if and only if the complex line subbundle  $\ell_2 := \pi^{\bot}_{\ell_1} (\Pz \ell_1) \subset \Sigma \times \C^5$ is isotropic. This is clearly equivalent to the orthogonality of $\overline{\ell_2}$ and $ \ell_2$. It is easily seen that a conformal map $f$ is isotropic if and only if $\la f_{zz}, f_{zz} \ra^c  \equiv 0$ on any complex coordinate $z$. On the other hand  a conformal map $f$ is called superconformal if  $\la f_{zz}, f_{zz} \ra^c   \not\equiv 0$. \\
If $f$ is isotropic then  $\ell_0, \ell_1, \ell_2$ are mutually orthogonal line subbundles of the trivial bundle $\Sigma \times \bb{C}^5$, whereas
 $ \ell_1 \oplus \ell_2$ is a maximal isotropic subbundle. The relevant geometric information of an  isotropic immersion $f : \Sigma \to \bb{S}^4$ is thus  contained in the ordered $3$-uple $(\ell_0, \ell_1, \ell_2)$. This motivates introducing the manifold $\Fl$ consisting of ordered $3$-uples $(X_0, X_1, X_2)$ of mutually orthogonal complex lines in $\C^{5}$ where  $X_0$ is the complexification of a real line in $\R^5$ and  $ X_1, X_2$ span a maximal complex isotropic subspace of $\C^5$.
 Given  $l = (X_0, X_1, X_2) \in \Fl$,  choose an ordered orthonormal basis  $\{u_0,   u_1,  \dots u_4 \}$ of $\R^{5}$ so that
\begin{equation}\label{projectonF} 
 X_0 = \C u_0, \quad  X_1 = \C z_1, \quad X_2 =  \C z_2, 
 \end{equation} 
   where  the complex unit vectors $z_1, z_2$ are given by
\begin{equation}\label{unitzj}   
z_j = \frac{u_{2j-1}- i u_{2j}}{\sqrt{2}},\, j=1,2
\end{equation} 
   
Define the projection map  $\pi : \Fl \to \bb{S}^4$   by demanding that the orthonormal basis $\{\pi(l),   u_1,  \dots u_4 \}$ be positively oriented. According to this definition $\pi(l) = \pm u_0 \in \bb{S}^4$, depending on the orientation of the chosen orthonormal basis.      
 Moreover every element $l= (X_0, X_1, X_2) \in \Fl$ is uniquely determined  by its projection $\pi(l) $ and the  complex isotropic line $X_1$.
 This establishes a diffeomorphism between the Grassmann bundle $G_2(T\bb{S}^4)$ and $\Fl$. In fact let $p \in \bb{S}^4$ and $V$ an oriented  $2$-plane in $T_p \bb{S}^4$. Choose an oriented orhtonormal base $\{ e_1, e_2 \}$  in $V$ and set  $X_0 := \C p$, $X_1 := \C(e_1-ie_2)$. Thus the application   sending $(p,V)$ to the (uniquely determined)  element $\ell =(X_0, X_1, X_2) \in \Fl$ is the desired   diffeomorphism.  \\         
Note that the projection  $\pi : \Fl \to \bb{S}^4$ is a fiber bundle, where the fiber $\Fl_p$ over a point $p \in \bb{S}^4$ is the totality of complex isotropic  lines in $T^{\C}_p \bb{S}^4$. Thus $\Fl_p$ identifies with    the Grassmann manifold $G_2 (\R^4) = SO_4 / SO_2 \times SO_2 \equiv \bb{S}^2 \times \bb{S}^2$,  a complex manifold. \\
 If $f: \Sigma \to \bb{S}^4$ is a conformal isotropic map we define its Gauss map $\widehat{f}: \Sigma \to \Fl$    by 
\begin{equation}\label{Gaussliftisotropic}
  \widehat{f} = (\ell_0, \ell_1, \ell_2 )
  \end{equation}
  It is easily verified that with this definition $\pi \circ \widehat{f} =f $ holds.\\  
On the other hand if $f$ is superconformal then $\overline{\ell}_2$ and $\ell_2$ are no longer  orthogonal. In this case we define $\widehat{f}: \Sigma \to \Fl$  by the condition   $\pi \circ \widehat{f} =f$. That is, 
there is a uniquely defined isotropic complex line subbundle $X_2 \subset \Sigma \times \C^5$  satisfying  $\pi (\ell_0, \ell_1, X_2)= f$. Note that $\ell_2$ and $\overline{\ell}_2$ are  contained in $X_2 \oplus \overline{X}_2$ which coincides with  the orthogonal complement bundle of $\overline{\ell}_1 \oplus \ell_0 \oplus \ell_1$. 
Define the Gauss map of a superconformal immersion $f$  by 
\begin{equation}\label{Gaussliftsupconf}
 \widehat{f} = (\ell_0, \ell_1, X_2)
 \end{equation}

Let $(X_0, X_1, X_2) \in \Fl $ with $X_0 =\C u_0$, and  $u_0 \in \R^5$ and $X_j = \C z_j, \, j=1,2$ where  the unit complex vectors $z_1, z_2 \in \C^5$ are given by~\eqref{unitzj}.  If $g \in SO_5$ then  
$$g.z_j = \frac{g. u_{2j-1}- i g. u_{2j}}{\sqrt{2}}, \quad j=1,2 $$
are also complex unit vectors $\in \C^5$. Therefore the following defines  a transitive action of the Lie group  $SO_5 $  on $\Fl$  
   $$         g. (X_0, X_1, X_2) =   (\C g.z_0,  \C g.z_1, \C g.z_2). $$ 
Hence  $\Fl \equiv  SO_5 /SO_2 \times SO_2 $, where the isotropy group $SO_2 \times SO_2$ is the  stabilizer of the basepoint  
 $$o=(\C e_0,  \C (\frac{e_1 - i e_2}{\sqrt{2}}),  \C (\frac{e_3 - i e_4}{\sqrt{2}})) \in \Fl.$$

\vskip .5cm

Let    $X \in \goth{g} = \goth{so}_5$, so that   $g_t= \exp (tX) \in SO_5, \, \forall t$. Given  
$l= (\C u_0,  \C z_1, \C z_2) \in \Fl$, with $z_j$ give by~\eqref{unitzj}. 
  Then $t \mapsto g_t . l $ defines a curve in $\Fl$ through $l$. We compute its derivative at $t=0$, 
\begin{equation*}
\begin{array}{cc}
\\
   \frac{d}{dt}|_{t=0} \exp (tX). (\C u_0,  \C z_1, \C z_2) = \\
   \\
   \frac{d}{dt}|_{t=0} (\C \exp (tX) u_0,  \C \exp (tX) z_1, \C \exp (tX) z_2)= \\
   \\
         (\C X u_0,  \C X(u_1 - i u_2), \C X(u_3 - i u_4)) \in T_{\ell } \Fl.\\
         \\
\end{array}
\end{equation*}   
Since $X+X^T=0$ it follows that 
$$Xu_0=0, \quad  \la Xu_k , u_k\ra = 0, \quad k=1,2,3,4. $$
On the other hand taking derivative at $t=0$ of the identities 
$$ \left\la   \exp (tX) \left( \frac{u_{2k-1}-i u_{2k}}{\sqrt{2}} \right),   \exp (tX) \left( \frac{u_{2k-1}-i u_{2k}}{\sqrt{2}} \right) \right \ra  = 1, \quad  \forall t \in \R, k=1,2,$$   
yields
$$    \la X u_{2k-1} , u_{2k} \ra = \la X u_{2k},   u_{2k-1} \ra =0, \quad  k=1,2.$$

 In particular any  tangent vector at the basepoint $o \in \Fl$ is an ordered  $3$-tuple of complex lines of the form 
 $$           ( \C Xe_0,  \C X(e_1- i e_2) , \C X(e_3-ie_4 )), $$

 where  $X\in \goth{so}_5$ varies in  the subspace $\goth{p} \subset \goth{so}_5 $ of matrices  of the form
   \begin{equation}\label{matrixinP}
         \begin{pmatrix}
                0 & a & b & c & d \\
               -a & 0 & 0 & e & f \\
               -b & 0 & 0 & g & h \\
               -c & -e& -g& 0 & 0 \\
               -d & -f& -h& 0 & 0\\
               \end{pmatrix}
               \end{equation} 
 In this way  $T_o \Fl$ identifies with $\goth{p}$.               
 \vskip .5 cm  
 
If  $g \in SO_5$ let $\tau_g : \Fl \to \Fl$ be the isometry sending $g'.o $ to $gg'.o $. It easily follows  that the projection $\pi : \Fl \to \bb{S}^4$ satisfies
 \begin{equation}\label{equivariantPI} 
                     \pi \circ \tau_g  = \tau'_g \circ \pi, \quad \forall   g \in SO_5,
 \end{equation}
 where $\tau'_g$ is the isometry of $\bb{S}^4$ induced by $g \in SO_5$. \\
 
\vskip .5cm

\subsection{The Maurer-Cartan form of $\Fl$.}  Decompose $\goth{so}_5 =\goth{g}=  \goth{k} \oplus \goth{p}$, where $\goth{k}= \goth{so}_2 \oplus \goth{so}_2$ is the Lie algebra of the maximal compact 2-torus $T^2 = SO_2 \times SO_2$ sitting within $SO_5$ according to the inclusion
\begin{equation}\label{inclusionKSO5}
 SO_2 \times SO_2 \ni (A,B) \mapsto diag (1,A,B) \in SO_5.
\end{equation}
 Let  $\goth{p} \equiv T_o \Fl$ be the set of real skew-symmetric matrices~\eqref{matrixinP}.
  For $A,B \in \goth{p}$ define 
  \begin{equation}\label{normalmetric}
          \la A,B\ra = -\frac{1}{2}trace (AB)
          \end{equation}         
Then $\la . , . \ra$ is an $Ad(T^2)$-invariant inner product on $\goth{p}$ which determines an $SO_5$-invariant metric on $\Fl$, the so-called normal metric. \\
 
 The  geometry of the flag manifold $(\Fl, \la \, , \, \ra) $ may be studied  with  aid of the so-called  Maurer-Cartan one form $\beta : T \Fl \to \goth{so}_5$  introduced by Burstall and Rawnsley in~\cite{burstall-rawnsley}  which we describe below.  \\

Recall   the
reductive decomposition $ \goth{g}=\goth{so}_5= \goth{k} \oplus \goth{p}$. Consider the surjective application  $\xi_o : \goth{g}  \ni X \mapsto  \frac{d}{dt}|_{t=0} \exp(tX).o \in T_o \Fl$. Thus  $\xi_o$ has  kernel  $\goth{k}$ and  restricts to an isomorphism $\goth{p} \to  T_o \Fl$. 
Form the associated vector bundle $ [\goth{p}] := SO_5 \times_K \goth{p}$, then  the map 
$$    [(g, X)] \mapsto \frac{d}{dt}|_{t=0} \exp(tAd(g)X).\,x =  d \tau_g (\frac{d}{dt}|_{t=0} \exp(tX).o ), \quad x= g.o,$$
 establishes an isomorphism of the associated bundle $[\goth{p}]$ and the tangent bundle $T \Fl$, where $\tau_g$ is the isometry of $\Fl$ sending $g'.o$ to $gg'.o$.\\ 
 
Since $\goth{p}$ is an $Ad(K)$-invariant subspace of  $\goth{g}$,  one has the inclusion $[\goth{p}] \subset [\goth{g}]:= \Fl \times \goth{g}$,  given by $ [\goth{p}] \ni [(g,X)] \mapsto (g.o, Ad(g)X) \in [\goth{g} ]$. 
Note that the fiber of $[\goth{p}] \to \Fl$ over the point $g.o$ identifies with 
$ \{ g.o\} \times Ad(g)\goth{p} \subset [\goth{g}]$.  Hence   
 there exists an identification of  $T \Fl$  with a subbundle of the trivial bundle $[\goth{g}]$. This inclusion   may be viewed as an $\goth{g} $-valued one-form  on $\Fl$ which will be denoted by $\beta$.   
Every   $X \in \goth{g} $ determines a flow on $\Fl$ defined by $\varphi_t (x) = \exp(tX).x$, which  is an isometry of $\Fl$ for any $t \in \R$. The vector field of the flow is then a Killing field denoted by $X^*$ which is  defined  by 
$$ X^*_x =   \frac{d}{dt}|_{t=0}
\exp(tX). x, \forall x \in \Fl.    
$$   
It is not difficult to show that 
$$    \beta_x(X^*) = Ad(g)[Ad(g^{-1})(X)]_{\goth{p}}, \quad \forall X \in \goth{p}, \quad x=g.o \in \Fl$$
 In particular at $o \in \Fl$ we have $ \beta(X^*_o)= X$ for any $X \in \goth{p}$.  From this formula it follows the equivariance of $\beta$ which is expressed by 
 \begin{equation} \label{equivariantbeta}
    \beta \circ d \tau_g = Ad(g) \beta, \quad \forall g \in SO_5.
    \end{equation}
   Equivalently, note that $\xi$ in turn satisfies 
  $$d \tau_g \circ \xi_o  = \xi_{g.o} \circ Ad(g), \forall g \in SO_5.$$
    
  \vskip .5cm  
For   $x = g.o \in \Fl$  the application  
  $\xi_x :
 \goth{g}  \to T_{x} \Fl$ such that 
$
        X \overset{\xi_x}{\longmapsto} X^*_x$, 
 maps $\goth{g}$ onto $ T_{x} \Fl$,    and  
  restricts to an isomorphism $ Ad(g)(\goth{p})   \to  T_x \Fl$ whose  inverse  conicides with 
  $\beta_{x}$. More details and properties of the one form $\beta$ and  proofs can be found in~\cite{burstall-rawnsley}.\\
   
\begin{remark}\label{metriconp}
We fix the bundle metric on $[\goth{p}]$ for  which  $\beta_{F.o} : T_{F.o} \Fl \to Ad(F)(\goth{p})$ an isometry for every $F \in SO_5$.  
\end{remark}

\section{ Normal harmonic maps into $\Fl$}\label{normalharmonicmaps}  

Here we study the harmonic map equation for smooth maps 
   from a Riemann surface $\Sigma$ into the flag manifold  $\Fl=SO_5/T^2$, on which  we have fixed  the normal metric $\la .,. \ra$ defined by~\eqref{normalmetric}. Our approach is based on~\cite{black} (see also~\cite{higaki}) and  is also valid  for a wider class of  naturally  reductive  homogeneous spaces.
Assume that $ \Sigma$ is compact and consider the energy of $\phi$ by
\begin{equation}\label{energyonomega}
E (\phi) = \frac{1}{2}\int_{\Sigma} \| d \phi \|^2 dA_g,
\end{equation}
   where $dA_g$ is the area form  on $\Sigma$ determined by a  conformal metric $g$, and  $\| d \phi \|^2$ is the Hilbert-Schmidt norm of $d \phi$ defined by  $\| d \phi \|^2 = \sum_i \la d \phi (e_i), d \phi (e_i) \ra$ for any orthonormal frame $\{ e_i \}$ on $\Sigma$.
 By definition $\phi$ is harmonic if it is an extreme of the energy functional $ \phi \mapsto     E (\phi)$.    \\
 
Denote  by  $\nabla$ the Levi-Civita connection on $\Fl$ determined by the normal metric  $\la.,.\ra$ and  by $\nabla^{\phi}$ the induced connection on the pull-back bundle $\phi^{*} T \Fl \to \Sigma$.  As is well known, by the  the first variation formula of the energy~\cite{eells-lemaire} it follows that $\phi$ is harmonic if and only if  its tension vanishes, $$trace (\nabla d \phi) =0.$$ 
Being $\Sigma$ a  Riemann surface then     $\phi$ is harmonic  if and only if  on every  local complex coordinate $z$ on $\Sigma$ the following equation holds
\begin{equation} \label{generalharmoniceq}
        \nabla^{\phi}_{\Bz} d \phi(\Pz) = 0,        
\end{equation}
where the left hand of this equation is just  a non-zero multiple of the  tension field of  $\phi$~\cite{eells-lemaire}.   \\
 
  For our purposes we need to reformulate  equation~\eqref{generalharmoniceq}     in terms of  the Maurer-Cartan form or Moment map $\beta$ of $\Fl$, see~\cite{black},~\cite{burstall-rawnsley}. \\   
  
   Let  $D$ be the {\it canonical connection of second kind} i.e. the affine connection  on $\Fl$ determined by the condition that  the $D$-parallel transport along the curve  $t \to \exp(tX).x$ is realized  by $ d \exp(tX) $.  At the basepoint $o \in \Fl$ we have 
   \begin{equation} \label{canonicconnectionD}     D_{X^*} Y^*(o) = \frac{d}{dt}|_{t=0}  d \exp(-t X) Y^* = [X^*, Y^*](o) = -[X,Y]_{\goth{p}}. 
   \end{equation}
   
   Note also that $D$ is determined by  the condition  $(D_{X^*}  X^*)_o =0, \forall X \in \goth{p}$.  Since $\nabla \neq D$ and  $D$ is metric i.e. $D \la \,, \, \ra =0$, it follows that  $D$  has non-vanishing torsion.  From\eqref{canonicconnectionD} we obtain   
  \begin{equation}\label{torsionD}
  T_o^D(X^*,Y^*) = - [X,Y]_{\goth{p}}, \quad X,Y \in \goth{p}.
  \end{equation}
 
  The following formula due to F. E.  Burstall and J. Rawnsley~\cite{burstall-rawnsley} allows  to compute  $D$ in terms of $\beta $ and the Lie algebra structure of $\goth{g}= \goth{so}_5$, 
\vskip .5 cm  
  \begin{lemma}~\cite{burstall-rawnsley}   
  \begin{equation} \label{betaformula}     \beta(D_X Y) = X \beta(Y) - [ \beta(X), \beta(Y)], \quad X,Y \in \mathcal{X} (\Fl). 
  \end{equation}
  \end{lemma}
\bigskip

Let  us now compute the Levi-Civita connection $\nabla$ of  the normal metric on $\Fl$. Since $P: SO_5 \to \Fl$ is a riemannian submersion, $P$ sends $\goth{p}$-horizontal geodesics in $SO_5$ onto $\nabla$-geodesics in $\Fl$. Hence $\nabla_{X^*}X^* =0$ for every $X \in \goth{p}$ which implies  $\nabla_{X^*}Y^* + \nabla_{Y^*}X^* =0$, for any $X,Y \in \goth{p}$. Since $\nabla$  is 
torsionless we deduce   
\begin{equation}\label{nablanormal}
 \nabla_{X^*}Y^*= \frac{1}{2} [X^*, Y^*], \forall X,Y \in \goth{p}.
 \end{equation}
 Hence at $o \in \Fl$ we get
 \begin{equation}\label{nablanormalato}
 (\nabla_{X^*}Y^*)_{o}= -\frac{1}{2} [X,Y]_{\goth{p}}, \,\, \forall X,Y \in \goth{p}.
 \end{equation}     
   We are now ready to obtain the following formula for the Levi-Civita connection $\nabla$ on $\Fl$ in terms of $\beta$  (for an equivalent formula see~\cite{ohnita-udagawa})
\begin{lemma}
   \begin{equation}\label{formulaLeviCivita}
   \beta (\nabla_X Y) =  X \beta(Y) - [ \beta(X), \beta(Y)] + \frac{1}{2} \pi_{\goth{p}} ([ \beta(X), \beta(Y)]), \quad X,Y \in \mathcal{X}(\Fl),
   \end{equation}
   where $\pi_{\goth{p}} : \Fl \times \goth{so}(5) \to \Fl \times_K \goth{p} =[\goth{p}] \equiv T \Fl$ is the projection onto the  the tangent bundle of $\Fl$.
   \end{lemma}
    \noindent {\bf Proof.} Let $X^*, Y^* $ be fundamental (Killing) vector fields on $\Fl$ determined by $X,Y \in \goth{p}$. \\
From de definition of $\beta$,  \eqref{canonicconnectionD} and \eqref{nablanormalato} we have 
    $$    
   \beta ((\nabla_{X^*} Y^*)_{o}) - \beta ((D_{X^*} Y^*)_{o}) = -\frac{1}{2} [X,Y]_{\goth{p}}+ [X,Y]_{\goth{p}}=
   \frac{1}{2} [X,Y]_{\goth{p}}, \,\,  \forall X,Y \in \goth{p}. $$ 
On the other hand the difference tensor $\nabla -D $ is  $SO_5$-invariant, and so is  $\beta(\nabla-D) = \beta(\nabla) - \beta(D)$ by formula ~\eqref{betaformula}. Hence it is determined by its value at the point $o \in \Fl$. For arbitrary vector fields on $\Fl$ formula~\eqref{formulaLeviCivita} follows from a straightforward  calculation. \hfill $\square$ \\

Define the $D$-fundamental form of a smooth map $\phi :\Sigma \to \Fl$ by 
    \begin{equation} \label{D2ndfundform}
               D d\phi (U,V) =  D^{\phi}_{U} d \phi(V) - d \phi (\nabla^{\Sigma}_U V), \quad  U,V \in \mathcal{X}(\Sigma),
    \end{equation}
   in which $D^{\phi}$ is the connection on the pullback bundle $\phi^{*} T \Fl \to \Sigma$ determined by $D$  and $\nabla^{\Sigma}$ is the Levi-Civita connection on $M$ determined by a conformal metric.  The map $\phi :\Sigma \to \Fl$ is called affine or $D$-harmonic if and only if   $trace (D d\phi) =0$ holds, or equivalently 
   $$D^{\phi}_{\Bz} d \phi(\Pz)=0.$$ 
 \begin{lemma}\label{harmonicmaplemma}     $\phi :\Sigma \to \Fl $ is normal-harmonic if and only if   it is $D$-harmonic. Hence $\phi$ is harmonic if and only if  
 \begin{equation} \label{Zharmonicmapequation}
 \bz (\phi^* \beta)' -
 [(\phi^* \beta)''  \wedge (\phi^* \beta)' ] =0,
 \end{equation}
 where $\phi^* \beta = (\phi^* \beta)' + (\phi^* \beta)''$ is the decomposition of the (complex) one form $\phi^* \beta$ into its $(1,0)$ and $(0,1)$ parts.   
 \end{lemma}
 \noindent{\bf Proof.}  Follows as consequence of the following  formula for the tension of $\phi$ 
  \begin{equation} \label{betatorsionphi}
  \begin{array}{cc}
 -\beta \left( \nabla^{\phi}_{\Bz} d \phi (\Pz) \right) dz \wedge d \bar{z} = \bz (\phi^* \beta)' -
 [(\phi^* \beta)''  \wedge (\phi^* \beta)' ].\\
 \end{array}
 \end{equation}
 To obtain formula \eqref{betatorsionphi}, first note   that~\eqref{formulaLeviCivita} implies $\beta (tr \nabla d \phi) =  \beta (tr D d \phi)$, which since $\Sigma$ is a Riemann surface is equivalent to 
 $$      \beta(\nabla^{\phi}_{\Bz} d \phi (\Pz) ) = \beta(D^{\phi}_{\Bz} d \phi (\Pz)).  $$     
 On the other hand  from formula~\eqref{betaformula} we obtain   
 \begin{equation} \label{betatorsionphi2}
 \beta(D^{\phi}_{\Bz} d \phi (\Pz)) = \Bz \beta d \phi(\Pz) - [ \beta d \phi(\Bz), \beta d \phi(\Pz) ],      
 \end{equation}   
 from which \eqref{betatorsionphi} follows. $\square$\\

    Now we want to express the harmonic map equation in terms of the one form $\alpha = F^{-1} dF = \alpha_{\goth{k}} + \alpha_{\goth{p}}$, where $F$ is  a frame of $\phi$. From the identity 
    $$ \phi^* \beta = Ad (F) \alpha_{\goth{p}} $$
     we get   $$(\phi^* \beta)' = Ad (F) \alpha'_{\goth{p}}, \quad (\phi^* \beta)'' = Ad (F) \alpha''_{\goth{p}}.$$ 

\vskip .5 cm
   Now we use the identity (see~\cite{burstall-pedit} pag. 241) 
   $$       d [Ad(F)] = Ad(F) \circ ad \alpha, $$
   to compute  
   $$   d (\phi^* \beta) =   d [Ad(F) \alpha_{\goth{p}}] =  Ad(F) \{ d \alpha_{\goth{p}} + [ \alpha \wedge \alpha_{\goth{p}}] \}.$$
 Thus 
    \begin{equation*}
 \begin{array}{cc}
      \bz (\phi^* \beta)' = \bz \{  Ad(F) \alpha'_{\goth{p}}    \} =  Ad(F) \{   \bz \alpha'_{\goth{p}} + [\alpha''_{\goth{k}} \wedge \alpha'_{\goth{p}} ] + [\alpha''_{\goth{p}} \wedge \alpha'_{\goth{p}} ]   \}.\\
      \end{array}
      \end{equation*}

Now since  
 \begin{equation*}
 \begin{array}{cc} 
  [ (\phi^* \beta)'' \wedge (\phi^* \beta)'] = Ad(F)   [\alpha''_{\goth{p}} \wedge \alpha'_{\goth{p}} ] \\
 \end{array}
 \end{equation*}  
    
    We conclude   
$$  \bz (\phi^* \beta)' - [ (\phi^* \beta)'' \wedge (\phi^* \beta)'] = Ad(F) \{\bz \alpha'_{\goth{p}} + [\alpha_{\goth{k}}'' \wedge \alpha'_{\goth{p}} ]     \}.
$$    
\vskip .5cm

Thus a direct consequence of Lemma~\ref{harmonicmaplemma}  is the following  
\begin{proposition}
Let $\phi : \Sigma \to \Fl$ be a smooth map.  Then $\phi$ is  harmonic  if and only if 
 the $\goth{so}_5$-valued one form $\alpha = F^{-1}dF$ satisfies 
\begin{equation}\label{harmeqalpha}
        \bz \alpha'_{\goth{p}} + [\alpha''_{\goth{k}} \wedge \alpha'_{\goth{p}} ] =0.
\end{equation}
\end{proposition}

\medskip
Here note that $ [\alpha''_{\goth{k}} \wedge \alpha'_{\goth{p}} ]= [\alpha_{\goth{k}} \wedge \alpha'_{\goth{p}} ]$. \\ 

Set  $\alpha = A dz + B d \bar{z}$, where $A= F^{-1}F_z, B = F^{-1} F_{\bar{z}}$, then  decomposing $A = A_{\goth{k}} + A_{\goth{p}}$ and  $B = B_{\goth{k}} + B_{\goth{p}}$
we have   $ \alpha'_{\goth{p}} = A_{\goth{p}} dz$ and  $\alpha''_{\goth{p}} = B_{\goth{p}} d \bar{z}$. Thus respect to a complex local coordinate $z$ equation~\eqref{harmeqalpha} becomes 
 \begin{equation} \label{auxbeta}
 \bz \alpha'_{\goth{p}} + [\alpha''_{\goth{k}} \wedge \alpha'_{\goth{p}} ] = 
 - ( \Bz A_{\goth{p}} + [ B_{\goth{k}} ,  A_{\goth{p}}] \,) \,   dz  \wedge d \bar{z} 
 \end{equation}

\begin{corollary}
 $\phi : \Sigma \to \Fl$ is harmonic if and only if  for every frame $F$ of $\phi$ and any complex coordinate $z$  the complex matrices $A_{\goth{p}},B_{\goth{k}}$  satisfy 
\begin{equation}\label{harm.map.eq.phi}
    \Bz  A_{\goth{p}} + [ B_{\goth{k}}, A_{\goth{p}}] =0.
\end{equation}
\end{corollary}
\vskip .5cm
\begin{remark} 
From~\eqref{auxbeta} and~\eqref{betatorsionphi} we obtain the following  formula for the tension of $\phi$, where $F$ is any frame of $\phi$.
\begin{equation}\label{tensionphi2}
\beta \left( \nabla^{\phi}_{\Bz} d \phi (\Pz) \right)  = Ad(F)\{  \Bz  A_{\goth{p}} + [ B_{\goth{k}}, A_{\goth{p}}]  \},
\end{equation} 
  \end{remark}

 \section{ A formula for the tension of the Gauss map}\label{tensionGauss}
Here we obtain a formula for  the tension  of the Gauss map  of a conformal  immersion $f: \Sigma \to \bb{S}^4$. As a by-product we caracterize  conformal immersions with normal-harmonic Gauss maps. \\

Let   $F=(F_0, F_1, F_2, N_1, N_2) \in SO_5$ be   a (local) frame in column notation of the Gauss map $\widehat{f}$ of a conformal immersion $f$ i.e. $\widehat{f} = F.o $ on an open subset $U \subseteq \Sigma$.  If $\Sigma $ is simply connected or contractible then there always exists a global frame
 $F: \Sigma \to SO_5$.  By the above considerations it follows that  $F_1-iF_2$ is a local section of $\ell_1$ and $N_1 - iN_2$ is a local section of $\ell_2$ (resp. $X_2$) if $f$ is isotropic (resp. superconformal). In either case the Gauss map  of $f$ is given locally by     
 \begin{equation}\label{frameGauss}
    \widehat{f} =( \C F_0 , \C (F_1 - i F_2), \C (N_1 - i N_2 ))=F.o
  \end{equation}
In particular $F$ is also an adapted frame of $f$.  
Since    
$$  F' = diag(1, R(\theta_1), R(\theta_2)). F, \quad   R(\theta_i) = \begin{pmatrix}
                     \cos \theta_i & \sin \theta_i \\
                     -\sin \theta_i & \cos \theta_i
                     \end{pmatrix} ,\quad i=1,2$$
 is also a frame of $\widehat{f}$ and also an adapted frame of $f$, we may assume (after possibly aplying a gauge) that~\eqref{fadapted} holds on a local complex chart $z$. \\  
  Set $\alpha = F^{-1}dF =  A dz + B d \bar{z}$ as before and decompose  $A  = A_{\goth{k}} + A_{\goth{p}} $ and  $B  = B_{\goth{k}} + B_{\goth{p}} $ according to the decomposition $\goth{g} = \goth{k} \oplus \goth{p} $. Define the following complex one forms on $\Sigma$,
 \begin{equation}\label{complexoneforms} 
 \begin{array}{ll}
  \alpha_{\goth{k}}= A_{\goth{k}} dz + B_{\goth{k}} d\bar{z}, & 
 \alpha_{\goth{p}}= A_{\goth{p}} dz + B_{\goth{p}} d\bar{z},\\
 \alpha'_{\goth{p}}=A_{\goth{p}} dz, & \alpha''_{\goth{p}}=B_{\goth{p}} d\bar{z},\\
 \alpha'_{\goth{k}}= A_{\goth{k}} dz, &  \alpha''_{\goth{k}}= B_{\goth{k}} d\bar{z},\\
 \end{array}
 \end{equation}   
where                                   
 \begin{equation}\label{matrApBp}
A_{\goth{p}} =
\begin{pmatrix}
0  &   -\frac{e^u}{\sqrt{2}}  &   i\frac{e^u}{\sqrt{2}}  &    0  & 0\\

\frac{e^u}{\sqrt{2}}& 0  &  0 &  - a_1   & -a_2\\

-i \frac{e^u}{\sqrt{2}} & 0 & 0 & -i  b_1  &  -i b_2  \\

0 & a_1  &  i b_1 &       0       &     0        \\

0 & a_2  &  i b_2 &    0     &        0              \\        
\end{pmatrix}, 
\, 
B_{\goth{p}} =
\begin{pmatrix}
0  &   -\frac{e^u}{\sqrt{2}}  &   -i\frac{e^u}{\sqrt{2}}  &    0  & 0\\

\frac{e^u}{\sqrt{2}}& 0  &  0 &  - \bar{a}_1   & -\bar{a}_2\\

i \frac{e^u}{\sqrt{2}} &  0 & 0 & i  \bar{b}_1  &  i \bar{b}_2  \\

0 & \bar{a}_1  &  -i \bar{b}_1 &       0       &    0        \\

0 & \bar{a}_2  &  -i \bar{b}_2 &   0     &        0              \\       
\end{pmatrix}
\end{equation}

\vskip 1 cm 

\begin{equation}\label{matrAkBk}
A_{\goth{k}} =
\begin{pmatrix}
  0 &                &               &     &  \\

  &       0         &  i u_{z} &     &  \\

  & -i u_{z}   &          0 &      &    \\

 &                 &    &    0    &   \sigma\\

 &   &   &    -\sigma     &        0      \\       
\end{pmatrix}, \,
B_{\goth{k}} =
\begin{pmatrix}
  0 &                &               &     &  \\

  &       0         &  -i u_{\bar{z}} &     &  \\

  & i u_{\bar{z}}   &          0 &      &    \\

 &                 &    &    0    &   \bar{\sigma}\\

 &   &   &    -\bar{\sigma}     &        0      \\       
\end{pmatrix}
\end{equation}       

Set
\begin{equation}\label{matrixM}
M:=\Bz  A_{\goth{p}} + [ B_{\goth{k}}, A_{\goth{p}}] = 
\begin{pmatrix}
0  &   0  &   0  &    0  & 0\\

0& 0  &  0 &  - A_1   & -A_2\\

0 & 0 & 0 & -i  B_1  &  -i B_2  \\

0 & A_1  &  i B_1 &       0       &     0        \\

0 & A_2  &  i B_2 &    0     &        0              \\        
\end{pmatrix},
\end{equation}

where  the coefficients $A_1,A_2, B_1, B_2$  are given respectively by

 \begin{equation}
\left \{ \begin{array}{ll}
 A_1 = \left( \Bz a_1 + u_{ \bar{z}} b_1 + a_2 \bar{\sigma} \right),\\
\\
A_2 = \left( \Bz a_2 + u_{ \bar{z}} b_2 - a_1 \bar{\sigma} \right),\\
\\
B_1 = \left(  \Bz b_1 + u_{ \bar{z}} a_1 + b_2 \bar{\sigma} \right),\\
\\
B_2 =  \left( \Bz b_2 + u_{ \bar{z}} a_2 - b_1 \bar{\sigma} \right).\\
\end{array} \right.
\end{equation}
\vskip .5cm

In terms  of the parameters  $u, h_i, \xi_i$, the above coefficients are given by,   

\begin{equation}
\left \{ 
\begin{array}{l}
 A_1=\frac{e^{-u}}{\sqrt{2}} [ \Bz \xi_1 + \xi_2 \bar{\sigma}] + \frac{e^{u}}{\sqrt{2}}[\Bz h_1 + h_2 \bar{\sigma}],\\
\\
 A_2= \frac{e^{-u}}{\sqrt{2}} [ \Bz \xi_2 - \xi_1 \bar{\sigma}] + \frac{e^{u}}{\sqrt{2}}[\Bz h_2 - h_1 \bar{\sigma}],\\
\\
 B_1= \frac{e^{-u}}{\sqrt{2}}  [ \Bz \xi_1 + \xi_2 \bar{\sigma}] - \frac{e^{u}}{\sqrt{2}}[\Bz h_1 + h_2 \bar{\sigma}],\\
\\
 B_2= \frac{ e^{-u}}{\sqrt{2}} [ \Bz \xi_2 - \xi_1 \bar{\sigma}] - \frac{e^{u}}{\sqrt{2}}[\Bz h_2 - h_1 \bar{\sigma}].\\
\end{array} \right.
\end{equation}
\vskip .5cm

Using   Codazzi's equations,
\begin{equation}\label{codazzi}
\left \{
\begin{array}{l}
         e^{u} ( \Pz h_1 + h_2 \sigma) = e^{-u}(\Bz \xi_1 + \xi_2 \bar{\sigma}),\\
       e^{u} ( \Pz h_2 - h_1 \sigma) =  e^{-u} (\Bz \xi_2 - \xi_1 \bar{\sigma}),\\
\end{array} \right.
\end{equation}
the above cofficients become 

\begin{equation}\label{coeffAiBi}
\left \{ 
\begin{array}{l}
 A_1=\frac{e^{u}}{\sqrt{2}} ( \Pz h_1 + \Bz h_1 + h_2 (\sigma  + \bar{\sigma})),\\
\\
A_2= \frac{e^{u}}{\sqrt{2}} ( \Pz h_2 +\Bz h_2 - h_1 (\sigma +  \bar{\sigma})),\\
\\
 B_1= \frac{e^{u}}{\sqrt{2}}  ( \Pz h_1 - \Bz h_1 + h_2 (\sigma  - \bar{\sigma})),\\
\\
 B_2= \frac{ e^{u}}{\sqrt{2}} ( \Pz h_2 - \Bz h_2 - h_1 ( \sigma  - \bar{\sigma})).\\
\end{array} \right.
\end{equation}

\begin{remark}\label{remarkcoeffAiBi}
From~\eqref{coeffAiBi} it follows that   $A_1, A_2$ are $\R$-valued functions while $B_1, B_2$ are $i \R$-valued. 
\end{remark}

Using formula~\eqref{tensionphi2} for the Gauss map we get

$$\beta \left( \nabla^{\widehat{f}}_{\Bz} d \widehat{f} (\Pz) \right)  = Ad(F) M, $$
where $M$ is given by~\eqref{matrixM} and $F$ is  a  frame of $\widehat{f}$, i.e. $\widehat{f} = F.o \in \Fl$. The tension of the Gauss map at a point $x\in \Sigma$ is given by 
$$ 
     \nabla^{\widehat{f}}_{\Bz} d \widehat{f}(\Pz) (x)  = (Ad(F(x))M(x))^*_{F(x).o} = F(x)M(x)F(x)^T. (F(x).o) \in T_{F(x).o} \Fl,           $$
   where $(Ad(F(x))M(x))^* $ denotes the fundamental or Killing vector field on $\Fl$ determined by $Ad(F(x))M(x) \in \goth{g}$.  \\
  Since  the base point was chosen to be $o= (\C e_0, \C(e_1-ie_2), \C(e_3-ie_4)) \in \Fl$,  a straightforward calculation shows (droping the point $x \in \Sigma$) that 
     $$
\nabla^{\widehat{f}}_{\Bz} d \widehat{f}(\Pz) = FM.o =    F. (\C M_0 , \C (M_1-iM_2), \C (M_3-iM_4)),
$$
where $M_j$ is the $j$-th column of the matrix $M$. In terms of the columns $F_0, F_1, F_2, N_1, N_2$ of $F$ we obtain  
\begin{equation}\label{tensionGaussLift}
\begin{array}{cc}
    \nabla^{\widehat{f}}_{\Bz} d \widehat{f}(\Pz)  = \\
    \\
    ( 0 , \C [(A_1 + B_1)N_1 + (A_2 + B_2) N_2], 
    \C[(-A_1 + iA_2)F_1 + (-iB_1+B_2)F_2 ] ),\\
    \end{array}
\end{equation}

Now  recall equation~\eqref{normalderivH} giving the normal  derivative of the mean curvature vector, \begin{equation*}
 \nabla^{\bot}_{\Pz} H = ( \Pz h_1 + h_2 \sigma) N_1 +   ( \Pz h_2 - h_1 \sigma) N_2.\\
\end{equation*}

From~\eqref{coeffAiBi} above, we obtain  
\begin{equation*}
\sqrt{2} e^{u} \, \nabla^{\bot}_{\Pz} H = (A_1 + B_1) N_1 + (A_2 + B_2)N_2. 
\end{equation*} 
Thus from~\eqref{tensionGaussLift} we arrive at the following formula for (a non-zero factor of) the tension of $\widehat{f}$,
\begin{equation}\label{tensionGaussLift2}
\nabla^{\widehat{f}}_{\Bz} d \widehat{f}(\Pz)  = (0, \C [\nabla^{\bot}_{\Pz} H], \C [\Psi]), 
\end{equation} 

where $\Psi = (-A_1 + iA_2)F_1 + (-iB_1+B_2)F_2$.\\

 From these identities we conclude that $\widehat{f}$ is harmonic if and only if $\nabla_{\Pz}^{\bot}H =0$ and $\Psi =0$.  
Both equations  $\nabla^{\bot}_{\Pz} H=0$ and $\Psi =0$ are equivalent, since by~\eqref{coeffAiBi} and  Remark~\eqref{remarkcoeffAiBi} each equation separatelly is equivalent to $ A_1 = A_2 = B_1 = B_2 =0$. 
We have thus completed the proof of the follwing theorem, 
\begin{theorem}\label{harmonicgausslift}
Let $f : \Sigma \to \bb{S}^4$ be a conformal immersion of a Riemann surface and  $\widehat{f}:\Sigma \to \Fl$ its Gauss map. Then 
$\widehat{f}$ is harmonic with respect to the normal metric on $\Fl$ if and only if the surface $f$ has parallel mean curvature vector,  $\nabla^{\bot}H=0$.
\end{theorem}

\section{Complete integrability }\label{integrability} 
 
 We first  turn our attention to the following  property of  the Gauss map, 
 \begin{lemma}\label{specialpropertylemma}
Let $\widehat{f}:\Sigma \to \Fl$ be the Gauss map of a conformal immersion $f$. Then for any frame $F$ of $\widehat{f}$ the complex one forms   $\alpha'_{\goth{p}} , \alpha''_{\goth{p}}$ defined by~\eqref{complexoneforms} satisfy the following condition 
 \begin{equation}\label{specialproperty}
  [\alpha'_{\goth{p}} \wedge \alpha''_{\goth{p}}]_{\goth{p}} =0.
 \end{equation}
 \end{lemma}

 \noindent {\bf Proof.} Recall that for $\goth{g}$-valued one forms  $\alpha, \beta$ on $\Sigma$  their wedge product is defined by  
 $$
 [ \alpha \wedge \beta] (X,Y) = [ \alpha(X), \beta(Y)] + [\beta(X) ,\alpha(Y)  ]
 $$
\vskip .5cm
   
Locally $  [\alpha'_{\goth{p}} \wedge \alpha''_{\goth{p}}]_{\goth{p}}=[A_{\goth{p}}, B_{\goth{p}}]_{\goth{p}} \, dz \wedge d\bar{z}  $, thus in order to prove~\eqref{specialproperty} it is enough to   prove   $[A_{\goth{p}}, B_{\goth{p}}]_{\goth{p}}=0$, for an arbitrary local complex coordinate $z$ on $\Sigma$.   
 Since   $ [A_{\goth{p}} , B_{\goth{p}} ]_{\goth{p}}$ is  a $5 \times 5$    skew-symmetric matrix, we need only check that  $8 = \dim \goth{p}$  entries vanish.  We  perform  a  direct computation using the structure equations~\eqref{structureeqs2} and the explicit form  of matrices $A_{\goth{p}}, B_{\goth{p}}$~\eqref{matrApBp}. \\
From  the structure equations~\eqref{structureeqs2} and taking~\eqref{aibi} into account we obtain  
$$\la  f_z , [\Bz N_i]^{T} \ra^c = \la  f_{\bar{z}} , [\Pz N_i]^T \ra^c = e^{2u} h_i, \quad i =1,2,$$
where $ [\Bz N_i]^{T}$, (resp.$ [\Pz N_i]^T $) denote projection of $\Bz N_i$ (resp. $\Pz N_i$) onto the tangent bundle of the immersed surface. Both equations together imply the vanishing of entries $(0,3)$ and $(0,4)$. \\
Denote by $\nabla^{\Sigma}$ the Levi-Civita connection on $\Sigma$ determined by the induced conformal metric $g= f^* \la .,.\ra$. 
Using  again the structure equations~\eqref{structureeqs2}  we obtain   
\begin{equation*}
\begin{array}{ll}
\la f_z , \Pz F_i - \nabla^{\Sigma}_{\Pz} F_i \ra^c =0, & i=1,2.\\
\\
\la \Pz F_i - \nabla^{\Sigma}_{\Pz} F_i , [\Bz N_i]^{ T} \ra^c = \la \Bz F_i - \nabla^{\Sigma}_{\Bz} F_i  , [\Pz N_j]^T \ra^c =0, & 1 \leq i,j \leq 2.\\
\end{array}
\end{equation*} 
    The first equation above implies  the vanishing of  entries $(0,1)$ and  $(0,2)$, while the second equation implies the vanishing of entries $(1,3)$,  $(1,4)$, $(2,3)$ and $(2,4)$.  
     This completes the proof of the Lemma. \hfill $\square$\\

Let $\widehat{f}: \Sigma \to \Fl $ be the Gauss map of $f$. Take a local frame $F$ of $\widehat{f}$ and consider the one form  $\alpha = F^{*} \Theta = F^{-1}dF$. If $\Sigma$ is simply connected then there always exists a global frame $F: \Sigma \to SO_5$ of any smooth map  into $\Fl$. According to the splitting   $\goth{g}=\goth{so}_5 = \goth{k} \oplus \goth{p}$   
  the Maurer-Cartan equation   
$  d \alpha + \frac{1}{2} [\alpha \wedge \alpha]=0$,   
   decomposes into the following equations 
\begin{equation}\label{decompalpha}
\left \{
\begin{array}{ll}     
    \bz \alpha'_{\goth{p}} + [ \alpha_{\goth{k}} \wedge \alpha'_{\goth{p}}] +
    \pz \alpha''_{\goth{p}} + [ \alpha_{\goth{k}} \wedge \alpha''_{\goth{p}}]+[ \alpha'_{\goth{p}} \wedge \alpha''_{\goth{p}}]_{\goth{p}}=0,\\
    d\alpha_{\goth{k}} + \frac{1}{2}[ \alpha_{\goth{k}} \wedge \alpha_{\goth{k}}]+
    [ \alpha'_{\goth{p}} \wedge \alpha''_{\goth{p}}]_{\goth{k}}=0.\\    
\end{array} \right.
\end{equation}

Further, assume that   $\widehat{f}$  is harmonic, hence it satisfies   
$$ 0= \bz \alpha'_{\goth{p}} + [\alpha_{\goth{k}} \wedge \alpha'_{\goth{p}} ].  $$        
Since  by Lemma~\ref{specialpropertylemma} condition   $[\alpha'_{\goth{p}} \wedge \alpha''_{\goth{p}}]_{\goth{p}}=0$ holds,  the pair of equations~\eqref{decompalpha} reduces to  
\begin{equation*}
\left\{
\begin{array}{ll}
(\textbf{a}) &\pz \alpha''_{\goth{p}} + [ \alpha_{\goth{k}} \wedge \alpha''_{\goth{p}}]=0,\\
(\textbf{b}) & d\alpha_{\goth{k}} + \frac{1}{2}[ \alpha_{\goth{k}} \wedge \alpha_{\goth{k}}]+
    [ \alpha'_{\goth{p}} \wedge \alpha''_{\goth{p}}]=0.\\
\end{array} \right.
\end{equation*}

If   $\lambda \in \bb{C}$ with $|\lambda| =1$   set 
\begin{equation}\label{alphalambdaoneform} 
\lambda. \alpha : = \alpha_{\lambda} = \lambda^{-1} \alpha'_{\goth{p}} + \alpha_{\goth{k}} + \lambda \alpha''_{\goth{p}}.
\end{equation}
Due to $\overline{A_{\goth{p}}}=B_{\goth{p}}$ and  $\overline{A_{\goth{k}}}=B_{\goth{k}}$, $\alpha_{\lambda}$ is $\goth{g}$-valued for every $\lambda \in \bb{S}^1$.  
Moreover  $\lambda . \alpha = \alpha_{\lambda}$  defines an action of $\bb{S}^1$ on $\goth{g}^{\C}$-valued $1$-forms which leaves invariant the solution set of equations (\textbf{a}) and (\textbf{b}) above. 
Comparing coefficients of $\lambda $  it follows that  equations $(\bf a)$ and $(\bf b)$ above hold for $\alpha$ if and only if $\alpha_{\lambda}$ satisfies       
   \begin{equation}\label{zcc}  
  d \alpha_{\lambda} + \frac{1}{2} [\alpha_{\lambda} \wedge \alpha_{\lambda}]=0, \forall \lambda \in \bb{S}^1.  \end{equation}
 
 Due to  $\overline{\alpha'_{\goth{p}}} = \alpha''_{\goth{p}} $  and $\overline{ \alpha_{\goth{k}} }= \alpha_{\goth{k}}$, the one forms $\alpha_{\lambda}$ satisfy 
 $$          \overline{\alpha_{\lambda}} = \alpha_{\lambda},  \quad  \forall \lambda \in \bb{S}^1,               $$
 thus  $\alpha_{\lambda}$ is $\goth{g}$-valued. It determines a connection $d + \alpha_{\lambda}$ on the trivial $SO_5$ bundle over $\Sigma$ which by~\eqref{zcc} is  automatically flat.  
 By this reason equation~\eqref{zcc} is  called {\it zero curvature condition} (ZCC)~\cite{burstall-pedit}. 
In this way the harmonic map equation for Gauss maps to $\Fl$  is encoded in a loop of zero curvature equations, a  manifestation of complete integrability~\cite{Guest}. We summarize our discussion in the following

  \begin{theorem}\label{comp.int.Gaussmap}
  The harmonic map equation for the Gauss maps $\widehat{f}: \Sigma \to \Fl$ of  a conformal immersion $f: \Sigma \to \bb{S}^4$  can be expressed as a  loop of zero-curvature equations.
\end{theorem}

Conversely, let us assume that $\Sigma$ is simply connected (otherwise transfer the whole situation to the universal covering surface $\tilde{\Sigma}$).  Let $\alpha_{\lambda}$ be the loop of $\goth{so}_5$-valued 1-forms~\eqref{alphalambdaoneform}.   
Fixed a  point $m_o \in \Sigma$  we integrate to  solve  
\begin{equation}\label{extendedframeeqs}
dF_{\lambda} = F_{\lambda} \alpha_{\lambda}, \quad F_{\lambda} (m_o)= Id \in SO_5.
\end{equation}
 The solution      
$F_{\lambda}=(F_0(\lambda), F_1(\lambda), F_2(\lambda), N_1(\lambda), N_2(\lambda) ) :  \Sigma \to SO_5$  is  called an {\it extended frame}~\cite{burstall-rawnsley} and satisfies   
\begin{equation}\label{lambdamat}
F^{-1}_{\lambda} (F_{\lambda})_z = \lambda^{-1} A_{\goth{p}}+ A_{\goth{k}} , \quad F^{-1}_{\lambda} (F_{\lambda})_{\bar{z}} = \lambda B_{\goth{p}} + B_{\goth{k}}, \quad \forall \lambda \in \bb{S}^1. 
\end{equation} 
Since   
$$ (\alpha_{\lambda})_{\goth{p}} = \lambda^{-1} \alpha'_{\goth{p}} + \lambda \alpha''_{\goth{p}} =(\alpha_{\lambda})'_{\goth{p}} +(\alpha_{\lambda})''_{\goth{p}}, \quad  (\alpha_{\lambda})_{\goth{k}}= \alpha_{\goth{k}},$$   
     the one form $\alpha_{\lambda}$ satisfies equations $(\textbf{a})$ and $(\textbf{b})$ for every $\lambda \in \bb{S}^1$.  Thus if     $P : SO_5 \to \Fl$ denotes the projection map  $P (g) = g.o$, then    $\widehat{f}_{\lambda} =  P \circ F_{\lambda} :\Sigma \to \Fl$ is harmonic  $\forall \lambda \in \bb{S}^1$.
  The family  $\{ \widehat{f}_{\lambda}, \lambda \in \bb{S}^1 \}$ is called {\it the associated family} of the harmonic Gauss map $\widehat{f}$ (see~\cite{burstall-pedit}).
  Note that   $\widehat{f}_{ \{ \lambda =1 \}} = \widehat{f} $, hence  each $ \widehat{f}_{\lambda} $ is a  deformation of $\widehat{f}$.

\subsection{One parameter families of parallel mean curvature immersions.}
From~\eqref{lambdamat} we derive the following set of equations which encode the dependence of $F_{\lambda}$ from the parameter $\lambda$,  

\begin{equation} \label{lambdafn}
\left \{
\begin{array}{ll}
 \Pz F_0(\lambda) = \lambda^{-1} \frac{e^{u}}{\sqrt{2}} \left[ (F_1(\lambda)-i F_2(\lambda)   \right],\\
 \\
   \Pz F_1(\lambda)=
-\lambda^{-1} \frac{e^{u}}{\sqrt{2}}.F_0(\lambda) -i u_z F_2(\lambda) + \lambda^{-1} a_1 N_1(\lambda) + \lambda^{-1} a_2 N_1(\lambda), \\
\\
 \Pz F_2(\lambda)= i \lambda^{-1} \frac{e^{u}}{\sqrt{2}}.F_0(\lambda)+ i  u_z F_1(\lambda) + i\lambda^{-1}  b_1N_1(\lambda) + i\lambda^{-1}  b_2 N_2(\lambda),  \\
 \\
 \Pz N_1(\lambda) = -\lambda^{-1} a_1 F_1(\lambda) - i \lambda^{-1} b_1 F_2(\lambda) - \sigma N_2(\lambda)\\
\\
\Pz N_2(\lambda) = - \lambda^{-1} a_2 F_1(\lambda) - i\lambda^{-1} b_2 F_2(\lambda) + \sigma N_1(\lambda), \\ 
 \end{array} \right.
\end{equation}
\vskip .3 cm
where $a_i, b_i$ are defined by~\eqref{aibi}. 
From the first equation above we obtain  
 $$\la (f_{\lambda})_z ,(f_{\lambda})_z \ra =\la (f_{\lambda})_z ,(f_{\lambda})_{\bar{z}} \ra^c  = e^{2u}, \forall \lambda \in \bb{S}^1.$$
 
Thus  $\{ f_{\lambda}, \lambda \in \bb{S}^1 \}$  is a family of conformal immersions of $\Sigma$ into $\bb{S}^4$, with a common conformal factor $u$. In particular   all $f_{\lambda}$ induce the same metric for every $\lambda \in \bb{S}^{1}$ and consequently the same Gaussian curvature function i .e. $K_{\lambda}=K, \, \forall \lambda \in \bb{S}^1$.\\
Let $K^{\bot}_{\lambda}$ denote the normal curvature of $f_{\lambda}$. Then from the $4$th and $5$th equations in~\eqref{lambdafn} we obtain 
$$ 
 \sigma^{\lambda} = \la N_2(\lambda) , N_1(\lambda) \ra = \sigma,$$   
 
 which implies  $K^{\bot}_{\lambda} =K^{\bot}$ for all $\lambda$. \\
  
Denote by   $H_{\lambda}$  the mean curvature vector of $f_{\lambda}$. Since $u$ is the common conformal parameter of all $f_{\lambda}$ we get  
$$ H_{\lambda}= e^{-2u}  (f_{\lambda})^{\bot}_{\bar{z} z}. $$

\noindent Decomposing    
$$      H_{\lambda} = h^{\lambda}_1 N_1(\lambda) + h^{\lambda}_2 N_2(\lambda), \quad h^{\lambda}_i  := \la H_{\lambda}, N_i(\lambda) \ra, i=1,2$$ 
 and using~\eqref{lambdafn} we obtain        
\begin{equation*}
\begin{array}{cc}
H_{\lambda} = e^{-2u} \la (f_{\lambda})_{\bar{z}z}, N_1(\lambda) \ra N_1(\lambda) +
e^{-2u} \la (f_{\lambda})_{\bar{z}z}, N_2(\lambda) \ra N_2(\lambda)=\\
\\ -e^{-2u} \la (f_{\lambda})_{z}, \Pz N_1(\lambda) \ra N_1(\lambda)-e^{-2u} \la (f_{\lambda})_{z}, \Pz N_2(\lambda) \ra N_2(\lambda)=\\
\\
 h_1 N_1(\lambda) + h_2 N_2(\lambda). 
\end{array}
\end{equation*} 
Thus   $h^{\lambda}_i = h_i,\,\, i=1,2 $ and  so  the dependence of  $H_{\lambda}$  on $\lambda$ is  only through   $N_i (\lambda)$. Note that this implies $ \| H_{\lambda} \| = \| H \|$,  $ \forall \lambda \in \bb{S}^1$. \\
Also, since each $\phi_{\lambda}$ is harmonic then by Theorem~\ref{harmonicgausslift} it follows that each $f_{\lambda}$ has parallel mean curvature i.e. $\nabla^{ \bot}_{\lambda} H_{\lambda} =0$,  
 where $\nabla^{\bot}_{\lambda}$ denotes the covariant derivative of the normal bundle of $f_{\lambda} : \Sigma \to \bb{S}^4$.\\

\noindent Let   
 $\xi^{\lambda}_i : = \la (f_{\lambda})_{zz}, N_i(\lambda) \ra^c,  \,\, i=1,2$, then  from~\eqref{lambdamat} and~\eqref{lambdafn} we get 
\begin{equation}\label{hopflambda}
      \xi^{\lambda}_i = -\la (f_{\lambda})_{zz}, N_i(\lambda) \ra^c= - \la (f_{\lambda})_{z}, \Pz N_i(\lambda) \ra^c = \lambda^{-2} \xi_i,\quad  i=1,2. 
\end{equation}

Thus if  $II^{\lambda}$ denotes the second fundamental form of $f_{\lambda}$ then
 $$II^{\lambda} (\Pz, \Pz) = \lambda^{-2} \xi_1 N_1(\lambda) + \lambda^{-2} \xi_2 N_2(\lambda)$$
We summarize  the conclusions of this discussion in the following

\begin{proposition}\label{deformations}
Let $f: \Sigma \to \bb{S}^4$ be a conformal immersion of a simply connected Riemann surface $\Sigma$. If the  Gauss map  $\widehat{f} : \Sigma \to \Fl$ is  normal-harmonic,  then $f$ is part of an $\bb{S}^1$-loop of conformal immersions $f_{\lambda} : \Sigma \to \bb{S}^4$ with parallel mean curvature satisfying the following properties:\\
i)  all $f_{\lambda}$ have the same common conformal factor hence they  induce on $\Sigma$ the same metric for every $\lambda \in \bb{S}^{1}$ and consequently  the same Gaussian curvature function. Also all  $f_{\lambda}$ have the same normal curvature function given by~\eqref{normalcurvature}. \\

Let $F_{\lambda}=(F_0(\lambda), F_1(\lambda), F_2(\lambda), N_1(\lambda), N_2(\lambda) ))$
 be a global extended frame of $\widehat{f}_{\lambda}$, then \\
ii) The mean curvature vector of $f_{\lambda}$ is given by 
$$      H_{\lambda}= h_1 N_1(\lambda) + h_2 N_2(\lambda),          $$
where $H=h_1 N_1 + N_2 h_2$ is the mean curvature vector of $f= f_{\lambda=1}$. \\
iii) The  second fundamental form of $f_{\lambda}$ satisfies 
\begin{equation}\label{2ndfundlambda}
 II^{\lambda} (\Pz, \Pz) = \lambda^{-2} \xi_1 N_1(\lambda) + \lambda^{-2} \xi_2 N_2(\lambda),
 \end{equation}
where   $II (\Pz, \Pz) = \xi_1 N_1 + \xi_2 N_2$ is the  $(2,0) $ part of the second fundamental form $II$  of $f$. 
\end{proposition}

\section{On the normal energy }\label{normalenergy}

Here we compute the energy of the Gauss map $\widehat{f} : \Sigma \to \Fl$  of a a conformal immersion    $f: \Sigma \to \bb{S}^4$. The energy density of a smooth map $\phi :\Sigma \to \Fl$ is defined by the Hilber-Schmidt norm 
$$\| d \phi \|^2 (p) = \la d \phi (e_1),d \phi (e_1) \ra + \la d \phi (e_2), d \phi(e_2)\ra,$$ 
 where $\{ e_1, e_2 \}$ is any orthonormal basis of $T_p \Sigma$. \\

 To compute $\| d \widehat{f} \|^2$ we use   the induced metric   $g=f^*\la\, , \, \ra$ which is conformal  and  locally given by  $g = 2 e^{2u} dz \otimes d \bar{z}$, where $u$ is the conformal factor. Thus    
 $$ \frac{1}{2} \| d \widehat{f} \|^2 =   e^{-2u} \la \widehat{f}_z, \widehat{f}_{\bar{z}} \ra^c.$$
 
 Let $F$ be any frame of  $\widehat{f}$ then we have the following   identity 
 $$\widehat{f}^* \beta =  Ad(F) \alpha_{\goth{p}},$$
  in which  $\alpha_{\goth{p}} = A_{\goth{p}} dz + B_{\goth{p}} d \bar{z}$ and   $A_{\goth{p}}$ is given by~\eqref{matrA}, and $B_{\goth{p}} = \overline{A_{\goth{p}}}$. Using this and   Remark~\ref{metriconp}   we obtain
 \begin{equation} \label{energydensityf}
 \begin{array}{cc}
 \\
 \frac{1}{2} \| d \widehat{f} \|^2 =  e^{-2u} \la \widehat{f}_z, \widehat{f}_{\bar{z}} \ra^c =
    e^{-2u} \la \widehat{f}^* \beta (\Pz),  \widehat{f}^*\beta(\Bz) \ra^c =  \\
    \\
   e^{-2u} \la A_{\goth{p}}, B_{\goth{p}} \, \ra^c  
   = - e^{-2u} \frac{1}{2} tr(A_{\goth{p}}.B_{\goth{p}} )=\\
   \\  
      e^{-2u}  (e^{2u} + |a_1|^2 + |a_2|^2 + |b_1|^2 + |b_2|^2 )=  \\
      \\
      1 + \| H \|^2 + e^{-4u} (|\xi_1|^2 + | \xi_2|^2).
   \end{array}   
 \end{equation}
Therefore
\begin{equation}\label{energyfhat} 
       E(\widehat{f}\, ) = \int_{\Sigma}  [ 1 + \| H \|^2 + e^{-4u} (|\xi_1|^2 + | \xi_2|^2)] d A_g.
\end{equation}

From Gauss equation we have 
\begin{equation} \label{gauss2}    e^{-4u} (|\xi_1|^2+ | \xi_2|^2 ) = 1+ \| H \|^2 - K \geq 0,   \end{equation}
where $K$ is the Gaussian  curvature of the induced metric $ g=f^* \la.,. \ra$ on $\Sigma$. Hence 
\begin{equation*} 
       E(\widehat{f}\, ) = \int_{\Sigma}  [2( 1 + \| H \|^2 ) - K ]  d A_g.
\end{equation*}

On the other hand recall that the Willmore energy of the immersion $f$
is defined  by
\begin{equation}\label{willmoreenergy}
    W(f) = \frac{1}{2 \pi} \int_{\Sigma}  [1+ \| H \|^2 - K] d A_g \geq 0.
    \end{equation}
From~\eqref{gauss2} it follows  that   $W$  measures   how far is  $f$ from satisfying the condition $II( \Pz, \Pz) = f^{\bot}_{zz} =0$.  From~\eqref{covderivH} every   surface satisfying $II( \Pz, \Pz) =0$  has parallel mean curvature vector.  
Recall that surfaces  satisfying  $f^{\bot}_{zz} =0$  are exactly those for which the  ellipse of curvature at any point $p \in \Sigma$  reduces to   $H(p) \in T^{\bot}_p \Sigma$. 
 \\       
Combining~\eqref{gauss2} and~\eqref{willmoreenergy} we obtain  the following formula  relating   both energies
\begin{equation}\label{WE}
     E( \, \widehat{f} \,) = 4 \pi W(f)  + \int_{\Sigma} K dA_g.
     \end{equation}
 By Gauss-Bonnet formula $\int_{\Sigma} K dA_g = 2 \pi \mathcal{X}(\Sigma)$, where  $\mathcal{X}(\Sigma)$ is the Euler characteristic, hence  a topological invariant of the surface $\Sigma$. In terms of the genus  $g (\Sigma)$ the Euler characteristic is given by  $\mathcal{X}(\Sigma) = 2-2g (\Sigma)$. Since $W(f) \geq 0 $ for every conformal immersion $f$, we obtain a  lower bound for the normal  energy of Gauss maps,
 \begin{equation}\label{infboundE}
 E( \, \widehat{f} \,) \geq 2 \pi (2-2g (\Sigma))
 \end{equation}
An interesting question is to find classes of  immersions for which equality is attained in~\eqref{infboundE}. For instance if $\Sigma$ is homeomorphic to a sphere $\bb{S}^2$ then inequality~\eqref{infboundE} becomes $E(\widehat{f})  \geq 4 \pi$.  Let $f: \bb{S}^2 \to \bb{S}^4$ be  any great sphere or equatorial inclusion of the round unit sphere $\bb{S}^2$, then it  is totally geodesic, hence minimal with $K=1$. Using formula~\eqref{energyfhat} above we obtain   $E(\widehat{f}) = 4 \pi$. \\

The situation for tori is different. For,  let $\Sigma$ be homeomorphic to a $2$-torus $\bb{S}^1 \times \bb{S}^1$ then~\eqref{infboundE} becomes $E(\widehat{f})  \geq 0$. Thus  $E(\widehat{f}) =0$ if and only if $W(f) =0 $. This is equivalent to the condition  $f_{zz}^{\bot} \equiv 0$ which in turn implies  that the Gauss equation reduces to $K=1 + \|H \|^2 $, and so  integrating we get $0=\int_{\Sigma}K dA \geq \int_{\Sigma} dA >0$.  This shows that there is no conformal immersion of a $2$-torus into $\bb{S}^4$ with  zero normal energy Gauss map. The search for a positive lower bound for the normal energy of Gauss maps of conformally immersed tori is an interesting problem. Recently A. Gouberman and K. Leschke~\cite{gouberman-leschke} constructed a $\bb{CP}^3$-family of Willmore tori conformally immersed in $\bb{S}^4$ with zero spectral genus and Willmore energy $W= 2 \pi^2 n$, with $n$ a positive integer.

\end{document}